\documentclass[final,3p,11pt, times]{article}
\usepackage{authblk}
\usepackage{amsmath}
\usepackage{amssymb,latexsym}
\usepackage[german,english]{babel}
\usepackage{url}
\usepackage{graphicx}
\usepackage{gastex}
\usepackage{longtable}
\usepackage{lscape}
\usepackage{tabularx}
\usepackage{multicol}
\usepackage{verbatim}
\usepackage{multirow}
\usepackage{graphicx}
\usepackage{epsfig,graphics,graphicx}
\usepackage{geometry}
\usepackage[usenames]{xcolor}
\usepackage{tikz}
\newcommand\undermat[2]{%
	\makebox[0pt][l]{$\smash{\underbrace{\phantom{%
					\begin{matrix}#2\end{matrix}}}_{\text{$#1$}}}$}#2}

\hbadness=\maxdimen

\newtheorem{rem}{{\bf Remark}}

\newtheorem{lem}{{\bf Lemma}}[subsection]
\newtheorem{thm}{{\bf Theorem}}[section]
\newtheorem{prop}{{\bf Proposition}}[section]
\newtheorem{cor}{{\bf Corollary}}[section]

\usepackage{chngcntr}
\counterwithout{equation}{section}
\begin{document}
	
	\title{On the Spectral properties of power graphs over certain groups} 
	
	\author[1]{Yogendra Singh}
	\author[1]{Anand Kumar Tiwari}
	\author[2]{Fawad Ali}
	
	\affil[1]{\small Department of Applied Science, Indian Institute of Information Technology, Allahabad 211015, India}

	\affil[2]{Institute of Numerical Sciences, Kohat University of Science and Technology, Kohat 26000, KPK, Pakistan}
	
	\maketitle
	
	\hrule
	
\begin{abstract}
The power graph $P(\Omega)$ of a group $\Omega$ is a graph with the vertex set $\Omega$ such that two distinct vertices form an edge if and only if one of them is an integral power of the other. In this article, we determine the power graph of the group $\mathcal{G} = \langle s,r \, : r^{2^kp} = s^2 = e,~ srs^{-1} = r^{2^{k-1}p-1}\rangle$. Further, we compute its characteristic polynomial for the adjacency, Laplacian, and signless Laplacian matrices associated with this power graph. In addition, we determine its spectrum, Laplacian spectrum, and Laplacian energy.   		
\end{abstract}

\smallskip

\textbf{Keywords:} Power graph, spectral radius, Laplacian spectrum.

\textbf{MSC(2010):} 05C25, 05C50.

\hrule

\section{Introduction}

All the graphs investigated in this study are simple, that is, finite undirected graphs having no multiple edges and loops.
In \cite{14}, Kelarev and Quinn presented the notion of directed power graphs $\overrightarrow{P}(\mathbb{S})$ of a semigroup $\mathbb{S}$.  Such a graph has vertex set $\mathbb{S}$ and there is an edge from the vertex $u$ to $v$ iff $u \neq v$ and $v = u^l$ for some $l \in \mathbb{N}$.
Subsequently, in \cite{1}, the authors demonstrated the undirected power graph $P(\mathbb{S})$ in which two different vertices are edge connected whenever one is a power of the other. They further determined that for every finite group $\Omega$, the corresponding power graph is complete iff $\Omega$ is a cyclic group of order $1$ or $q^{\ell}$, where $q$ is any prime and $\ell \in \mathbb{N}.$ The authors in \cite{3} established that if  $\Omega_1$ and $\Omega_2$ are two finite abelian groups such that $P(\Omega_1) \cong P(\Omega_2)$, then $\Omega_1 \cong \Omega_2$. In \cite{2}, it is shown that  $P(\Omega_1) \cong P(\Omega_2)$ iff $\overrightarrow{P}(\Omega_1) \cong \overrightarrow{P}(\Omega_2)$. Curtin and Pourgholi demonstrated in \cite{8} that all finite cyclic groups have the largest clique and the most edges in their power graphs. In addition, several experts have examined various notions of power graphs, for instance, see \cite{AR(2022), Ali2022, Rather1}.
The first review paper on power graphs was published in $2013$ \cite{survey0}, while the most recent review article was published in $2021$ \cite{kumar2021recent}.

In the existing literature, the study of graph spectra of the power graph is a fascinating concept in several branches of mathematics, including group theory, Lie algebra, and ring theory.  In fact, graph spectra have several applications in quantum chemistry. Over the last decade, there has been a lot of interest in calculating the spectrum of algebraic graphs, see \cite{Ali(2020), AL(2019), 6, ghorbani1, Rather2}.
The authors of \cite{6a, 15aa} calculated the spectrum and the Laplacian spectrum of power graphs of Mathieu groups, cyclic, and dihedral groups.
Motivated by their work, we extend and discover the adjacency, Laplacian and signless Laplacian spectra, and energy of the power graphs of the group $\mathcal{G}$.

In the literature review, it has been noted that there are still several gaps regarding the identification of various spectral characteristics of power graphs of the group $\mathcal{G}$. 
The obvious explanation for this is that neither the construction of power graph for finite groups nor the derivation of appropriate formulations of graph characteristics for large classes of groups is adequately described or possible.
In this paper, we attempt to investigate these problems.

The remaining sections are arranged as follows: Section \ref{Sec2} provides definitions, symbols, and findings that are required later. In Section \ref{Sec4}, spectral properties, such as the adjacency spectrum, Laplacian spectrum, spectral radius, and the energies of power graphs are discussed.
%The conclusion and future work are presented in Section \ref{Sec6}.

\section{Preliminaries} \label{Sec2}

This section covers several fundamental graph-theoretic concepts and well-known results that we use throughout the article. We refer to see \cite{cha(2006), West(2009)} for more details.

Suppose $\Gamma$ is a simple undirected graph. The vertex and edge set of $\Gamma$ is represented by $V(\Gamma)$ and $E(\Gamma)$.
Two vertices $v_{1}$ and $v_{2}$ are connected if there is an edge between them, and they are symbolized as $v_{1}\sim v_{2},$ otherwise $v_{1}\nsim v_{2}.$
The degree (valency) $deg(u_{1})$ of a vertex $u_{1}$ is the set of vertices connected to $u_{1}$ in $\Gamma$.

Let $\Gamma_{1}$ and $\Gamma_{2}$ be two connected graphs. 
Then $\Gamma_{1} + \Gamma_{2}$ is the sum of $\Gamma_{1}$ and $\Gamma_{2}$ whose vertex and edge sets are $ V(\Gamma_{1}) \cup V(\Gamma_{2})$ and $ E(\Gamma_{1}) \cup E(\Gamma_{2}) \cup \{ c \sim d : c \in V(\Gamma_{1}),~ d \in V(\Gamma_{2})\}$, respectively. We denote $n\Gamma$ as the disjoint union of $n$ copies of $\Gamma$. If every pair of vertices of a graph $\Gamma$ on $n$ vertices is connected by an edge, then the graph is called complete and is denoted by $K_{n}$.

The adjacency matrix of $\Gamma$ is an $n \times n$ matrix of the form $A(\Gamma) = (z_{ij})$, where $z_{ij} = 1$ if $i$ is edge connected to $j$, and $0$ otherwise. 
The polynomial $\Delta(\Gamma, ~x) = \det(xI-A(\Gamma))$ is the characteristic polynomial of $\Gamma$. 
The eigenvalues of $\Gamma$ are basically the eigenvalues of $A(\Gamma)$ and are indicated by $\lambda_{t}(\Gamma)$, where $t = 1, 2, \ldots, n$.
The matrices $L(\Gamma)=D(\Gamma)-A(\Gamma)$ and $Q(\Gamma)=D(\Gamma)+A(\Gamma)$ are the Laplacian and the signless Laplacian matrices of $\Gamma$, respectively.
The (signless) Laplacian eigenvalues of $\Gamma$ are the eigenvalues of $(Q(\Gamma))$ $L(\Gamma)$ and are indicated by $(\delta_{t}(\Gamma))$ $\mu_{t}(\Gamma)$, where $t = 1, 2, \ldots, n$.
%Their multiset of eigenvalues is the signless Laplacian and the Laplacian spectrums of $\Gamma$, respectively.
Due to the symmetry of these matrices, all of its eigenvalues are real. 
Consequently, such eigenvalues may be listed as follows: $\lambda_{1}(\Gamma) \geq \lambda_{2}(\Gamma) \geq \cdots \geq \lambda_{n}(\Gamma)$,
%The Laplacian and the signless Laplacian are positive real semi-definite matrices, so their spectrum is real, and they are ordeblue as
$\mu_{1}(\Gamma) \geq \mu_{2}(\Gamma) \geq \dots \geq \mu_{n}(\Gamma)$ and $\delta_{1}(\Gamma)\geq \delta_{2}(\Gamma)\geq \dots \geq \delta_{n}(\Gamma)$, respectively the adjacency, Laplacian and the signless Laplacian eigenvalues of $\Gamma$.
Further details about these matrices can be seen in \cite{BH, Cvetkovic(2010)}.
The spectrum of $\Gamma$ is defined as the collection of all eigenvalues and their multiplicity. 
The spectral radius is the highest eigenvalue (sometimes termed the Perron-Frobenius eigenvalue) $\lambda_{1}(\Gamma)$ of $\Gamma$. If $m = |E(\Gamma)|$ and $n = |V(\Gamma)|$, then $LE(\Gamma)$ = $\sum_{i=1}^{l}|\mu_{i}-\frac{2m}{n}|$ is called Laplacian energy of $\Gamma.$ 

We denote the matrix $\mathcal{J}_n$ such that its diagonal elements are 0 and other elements are 1,  matrices $\mathcal{O}_n$ and $I_n$ as the zero and identity matrix of order $n$, respectively.

%\textcolor{blue}{Let $G$ be a group and $N$ be its normal subgroup. Then $G$ is called metacyclic if $N$ and $G/N$ are cyclic. H\"{o}lder [\cite{Zassenhaus(1999)}, Theorem 7.21] proved that a group $G$ is metacyclic if and only if $G$ is of the following form:}

%\textcolor{blue}{$$ G_{\alpha, \beta,\gamma, \delta} = \langle a, b \, : a^{\alpha} = b^{\delta}, ~ b^{\beta} = 1, aba^{-1}= b^{\gamma}\rangle,$$ $\text{ where } \alpha, \beta, \gamma, \delta \in \mathbb{N}, ~ gcd(\beta, \gamma) = 1,~ \beta ~|~(\gamma^\alpha-1), \text{ and } \beta ~|~ \delta(\gamma-1).$}

%\textcolor{blue}{The group $ G_{\alpha, \beta,\gamma} = \langle a, b \, : a^{\alpha} = b^{\beta} = 1, aba^{-1}= b^{\gamma}\rangle $ is called split metacyclic group.}

\smallskip

%\textcolor{blue}{Note that, $G_{2, n,n-1} \cong D_{2n},$  $G_{2,4n,2n-1} \cong SD_{8n},$ and $G_{2,2n,n,2n-1} \cong Q_{4n}.$ Thus, the dihedral, semidihedral, dicyclic groups are the special cases of $G_{\alpha,\beta,\gamma,\delta}.$ }

In this article, we focus our attention on the group $\mathcal{G}$, which is a non-abelian group of order $2(2^kp)$, where  $k \geq 2$ and $p \neq 2$ is a prime number. The presentation of $\mathcal{G}$ is given as: $$\mathcal{G} = \langle s,r \, : r^{2^kp} = s^2 = e,~ srs^{-1} = r^{2^{k-1}p-1}\rangle.$$ 
From the structure of $\mathcal{G}$ and its power graph, we have the following useful remark:

\begin{rem} \label{rem1}
We may divide the group $\mathcal{G}$ in the following manner:
$$H_0 = \{e, u = r^{2^{k-1}p}\},~
H_1 =  \{r^i : 1 \leq i \leq 2^kp-1\} \ \setminus \ \{u\},~
H_2 =  \{sr^{2t} {: 1 \leq t \leq 2^{k-1}p}\}, \text{ and }$$ $H_3 = \{sr^{2j+1} : 0 \leq j \leq 2^{k-1}p-1\}.$  As $\langle r\rangle \ \cong \ \mathbb{Z}_{{2^k}p}$, $P(\langle r\rangle) \cong P(\mathbb{Z}_{{2^k}p}).$
Since $o \left(sr^{2t} \right)=2$, every element of $H_2$ is adjacent with $e$ only and generates a $K_2$ graph. 
As there are $2^{k-1}p$ such elements, so we get $2^{k-1}p$ copies of $K_2$ graphs. 
Since $o(sr^{2j+1}) = 4$, so $\langle sr^{2j+1} \rangle$ has $\phi(4) =2$ generators, where $\phi$ is the Euler's function. 
Note that $\langle r^{2^{k-1}p} \rangle \ \leq \ \langle sr^{2j+1} \rangle$. This implies the generators of $\langle sr^{2j+1} \rangle$  are adjacent to the vertices of $\langle r^{2^{k-1}p} \rangle$. 
 This gives $2^{k-2}p$ number of $P\left( \langle r^{2^{k-1}p} \rangle \right) + K_2$ graphs, 
as there are $2^{k-1}p$ numbers divided between $2$ generators.
\end{rem}

Therefore, from the above structure, we have: 

\begin{thm} For the group $\mathcal{G}$, we have
$$P(\mathcal{G}) = P \left(\mathbb{Z}_{2^kp} \right) \cup 2^{k-1}p K_2 \cup 2^{k-2}p \left((P \langle r^{2^{k-1}p} \rangle) + K_2 \right).$$

%The power graph $P(\mathcal{G})$ of $\mathcal{G}$ is given below in Figure $1$. Thepower graph of $Z_{2^kp}$ is $K_{2^k p}$, shown as the shaded part of the figure.
\end{thm}

\section{Spectral properties} \label{Sec4} %Characteristic polynomial
\paragraph{}
In this part, we investigate the characteristic polynomials for the adjacency, Laplacian, and signless Laplacian matrices associated with the power graph $P(\mathcal{G})$ for the previously mentioned group $\mathcal{G}$. In addition, we determine its spectrum, Laplacian spectrum, and Laplacian energy.

\smallskip

Recall the following results that we use throughout the paper. 

\begin{thm} [\cite{BH}] \label{t6.1}
	If $\mathcal{A}, \mathcal{B}, \mathcal{C}$, and $\mathcal{D}$ are square matrices of the same order and $\mathcal{D}$ is invertible. Then \[\begin{vmatrix}
	\mathcal{A} & \mathcal{B}  \\
	\mathcal{C} & \mathcal{D} \end{vmatrix} = \begin{vmatrix}\mathcal{D}	\end{vmatrix} \begin{vmatrix} \mathcal{A} - \mathcal{B}\mathcal{D}^{-1}\mathcal{C} \end{vmatrix}. \]
\end{thm}

\begin{thm} [\cite{Cvetkovic(2010)}] \label{t6.2}
	Let $\Gamma$ be a  connected graph and $v \in V(\Gamma)$. Then 
	$$ \lambda_1(\Gamma-v) <  \lambda_1(\Gamma).$$
\end{thm}

\begin{thm} [\cite{Cvetkovic(2010)}] \label{t6.3}
	Let $\mathcal{A}, \mathcal{B}$ be two Hermitian matrices of the same order. Then
	
	$$ \lambda_{r}(\mathcal{A}+\mathcal{B}) \leq  \lambda_{s}(\mathcal{A}) + \lambda_{r-s+1}(\mathcal{B}), \text{ where } n \geq r \geq s \geq 1. $$
\end{thm}

The following propositions can be shown by using basic concepts of linear algebra 
\begin{prop}\label{l6.3}
	Let $\mathcal{A}$ and $\mathcal{C}$ be two $m \times m$ and $n \times n$ matrices, respectively, and let $\mathcal{B}$ and $\mathcal{D}$ be matrices of the appropriate sizes. Then \[\begin{vmatrix}
    \mathcal{A} & \mathcal{B}  \\
	\mathcal{O} & \mathcal{C} \end{vmatrix} = \begin{vmatrix}
	\mathcal{A} & \mathcal{O}  \\
	\mathcal{D} & \mathcal{C}\end{vmatrix} = \begin{vmatrix}\mathcal{A} 	\end{vmatrix} \begin{vmatrix} \mathcal{C} \end{vmatrix}. \]
\end{prop}

\begin{prop} \label{l6.4}
	Let $\mathcal{A}$ be a square matrix of order $n$ such that its diagonal elements are $x$ and other elements are $y$. Then $|\mathcal{A}| = (x+(n-1)y)(x-y)^{n-1}.$ 
\end{prop}

\begin{prop} \label{t6}
	A block diagonal matrix $\mathcal{A} = diag(\mathcal{D}_1, \mathcal{D}_2, \ldots, \mathcal{D}_n)$ is invertible if and only if each of its main-diagonal blocks is invertible, and in this case, $$\mathcal{A}^{-1} = diag(\mathcal{D}_1^{-1}, \mathcal{D}_2^{-1}, \ldots, \mathcal{D}_n^{-1}).$$
\end{prop}

\begin{lem} \label{l6}
	The characteristic polynomial of $P(\mathcal{G})$ is:
%\begin{align*}

\bigskip
	$\Delta(P(\mathcal{G}), x)  = x^{2^{k-1}p-1}(x+1)^{5.2^{k-2}p-3}(x-1)^{2^{k-2}p-1}\bigg[x^5 - (2^kp-2)x^4 - 5.2^{k-1}px^3 + \big(3.2^{2k-1}p^2-2^{k+2}p-2 \big)x^2+ \big(5.2^{2k-2}p^2-3.2^{k-1}p-1 \big)x - \big(2^{3k-2}p^3-2^{2k-2}p^2-2^kp \big)\bigg].$
%	\end{align*}
	
\end{lem}

\noindent{\textbf{Proof.}}  The adjacency matrix $A$ of $P(\mathcal{G})$ is 
$\begin{bmatrix}

\mathcal{B}_{2^kp} & \mathcal{C}_{2^kp}  \\ \\
\mathcal{C}^T_{2^kp} & \mathcal{D}_{2^kp}
\end{bmatrix}$, where

\smallskip
\hspace{-.65cm}
$\mathcal{C}_{2^kp} = \begin{array}{@{} c @{}}
\begin{bmatrix}
\begin{array}{ *{8}{c} }

\mathcal{E} & \mathcal{E} & \cdots & \mathcal{E} & \mathcal{F} & \mathcal{F} & \cdots & \mathcal{F} \\
\mathcal{O} & \mathcal{O} & \cdots & \mathcal{O} & \mathcal{O} & \mathcal{O} & \cdots & \mathcal{O} \\
\mathcal{O} & \mathcal{O} & \cdots & \mathcal{O} & \mathcal{O} & \mathcal{O} & \cdots & \mathcal{O} \\
\vdots & \vdots & \vdots & \vdots & \vdots & \vdots & \ddots  & \vdots \\
\mathcal{O} & \mathcal{O} & \cdots & \mathcal{O} & \mathcal{O} & \mathcal{O} & \cdots & \mathcal{O} \\

\undermat{2^{k-2}p}{\mathcal{O} & \mathcal{O} & \cdots & \mathcal{O}} & \undermat{2^{k-2}p}{\mathcal{O} & \mathcal{O} & \cdots & \mathcal{O}}  \\
\end{array}
\end{bmatrix} \\
\mathstrut
\end{array}$, $\mathcal{D}_{2^kp} = \begin{array}{@{} c @{}}
\begin{bmatrix}
\begin{array}{ *{9}{c} }

\mathcal{G} & \mathcal{O} & \mathcal{O} & \cdots & \mathcal{O} & \mathcal{O} & \cdots & \mathcal{O} \\
\mathcal{O} & \mathcal{G} & \mathcal{O} & \cdots & \mathcal{O} & \mathcal{O} & \cdots & \mathcal{O} \\
\mathcal{O} & \mathcal{O} & \mathcal{G} & \cdots & \mathcal{O} & \mathcal{O} & \cdots & \mathcal{O} \\
\vdots & \vdots & \vdots & \ddots & \vdots & \vdots & \ddots & \vdots \\
\mathcal{O} & \mathcal{O} & \mathcal{O} & \cdots & \mathcal{G} & \mathcal{O} & \cdots & \mathcal{O} \\
\mathcal{O} & \mathcal{O} & \mathcal{O} & \cdots & \mathcal{O} & \mathcal{O} & \cdots & \mathcal{O} \\
\vdots & \vdots & \vdots & \vdots & \vdots & \vdots & \ddots & \vdots \\
\undermat{2^{k-2}p}{\mathcal{O} & \mathcal{O} & \mathcal{O} & \cdots & \mathcal{O}} & \undermat{2^{k-2}p}{\mathcal{O} & \cdots & \mathcal{O}}  \\
\end{array}
\end{bmatrix} \\
\mathstrut
\end{array},$  \\ $$ \mathcal{B}_{2^kp} =
\begin{bmatrix}
0 & 1 & 1 & \cdots & 1 \\
1 & 0 & 1 & \cdots & 1 \\
1 & 1 & 0 & \cdots & 1 \\
\vdots & \vdots & \vdots &\ddots & \vdots \\
1 & 1 & 1 & \cdots & 0
\end{bmatrix},  \ \mathcal{E} =
\begin{bmatrix}
1 & 1 \\
1 & 1
\end{bmatrix}, \ \mathcal{F} =
\begin{bmatrix}
1 & 1 \\
0 & 0
\end{bmatrix},  \ \mathcal{G} =
\begin{bmatrix}
0 & 1 \\
1 & 0
\end{bmatrix},  \ \mathcal{O} =
\begin{bmatrix}
0 & 0 \\
0 & 0
\end{bmatrix}.$$
Thus, the characteristic polynomial of $P(\mathcal{G})$ is $\Delta(P(\mathcal{G}), x) =
\begin{vmatrix}
xI_{2^kp}-
\mathcal{B}_{2^kp} & -\mathcal{C}_{2^kp}  \\
-\mathcal{C}^T_{2^kp} & xI_{2^kp} - \mathcal{D}_{2^kp}
\end{vmatrix}$ 

Now, by Proposition \ref{t6}, we have

$$(xI_{2^kp} - \mathcal{D}_{2^kp})^{-1} = \begin{array}{@{} c @{}}

\begin{bmatrix}
\begin{array}{ *{9}{c} }

\frac{x}{x^2-1} & \frac{1}{x^2-1} & \cdots & 0 &  0 & 0 & \cdots & 0 & 0 \\
\frac{1}{x^2-1} & \frac{x}{x^2-1} & \cdots & 0 & 0 & 0 & \cdots & 0 & 0 \\
\vdots & \vdots & \ddots & \vdots & \vdots & \vdots & \ddots & \vdots & \vdots \\
0 & 0 & \cdots & \frac{x}{x^2-1} & \frac{1}{x^2-1} & 0 & \cdots & 0 & 0 \\
0 & 0 & \cdots & \frac{1}{x^2-1} & \frac{x}{x^2-1} & 0 & \cdots & 0 & 0 \\
0 & 0 & \cdots & 0 & 0 & \frac{1}{x} & \cdots & 0 & 0 \\
\vdots & \vdots & \ddots & \vdots &  \vdots & \vdots & \ddots & \vdots & \vdots \\
0 & 0 & \cdots & 0 & 0 & 0 & \cdots & \frac{1}{x} & 0 \\
0 & 0 & \cdots & 0 &  0 & \undermat{2^{k-1}p}{0 & \cdots & 0 & \frac{1}{x}}  \\
\end{array}
\end{bmatrix} \\
\mathstrut
\end{array}$$
So, by Theorem \ref{t6.1}, we have 
$$\Delta(P(\mathcal{G}), x)  = \begin{vmatrix}
xI_{2^kp}-\mathcal{D}_{2^kp}
\end{vmatrix}\begin{vmatrix}
xI_{2^kp}-\mathcal{B}_{2^kp}-\mathcal{C}_{2^kp}(xI_{2^kp} - \mathcal{D}_{2^kp})^{-1}\mathcal{C}^T_{2^{k}p}
\end{vmatrix}.$$

By using the fact that the determinant of a block diagonal matrix is the product of the determinant of its blocks, we have 

$$\begin{vmatrix}
xI_{2^kp}-\mathcal{D}_{2^kp}
\end{vmatrix} = (x^2-1)^{2^{k-2}p}x^{2^{k-1}p}.$$

Therefore, 
$$\Delta(P(\mathcal{G}), x)  =(x^2-1)^{2^{k-2}p}x^{2^{k-1}p} \begin{vmatrix}
xI_{2^kp}-\mathcal{B}_{2^kp}-\mathcal{C}_{2^kp}(xI_{2^kp} - \mathcal{D}_{2^kp})^{-1}\mathcal{C}^T_{2^{k}p}
\end{vmatrix}.$$

Now,
$$\mathcal{C}_{2^kp}(xI_{2^kp} - \mathcal{D}_{2^kp})^{-1}\mathcal{C}^T_{2^{k}p} = 
\begin{bmatrix}
\begin{array}{ *{9}{c} }

\frac{2^{k-1}p}{x} + \frac{2^{k-1}p}{x-1} & \frac{2^{k-1}p}{x-1} & 0 & \cdots & 0 & 0  \\ \\
\frac{2^{k-1}p}{x-1} & \frac{2^{k-1}p}{x-1} & 0 & \cdots & 0 & 0 \\
0 & 0 & 0 & \cdots & 0 & 0 \\
\vdots & \vdots & \vdots & \vdots & \vdots & \vdots \\
0 & 0 & 0 & \cdots & 0 & 0 \\
0 & 0 & \undermat{2^kp-2}{ 0 & \cdots & 0 & 0} \\
\end{array}
\end{bmatrix} \\
\mathstrut.$$

Therefore,

$$\begin{vmatrix}
xI_{2^kp}-\mathcal{B}_{2^kp}-\mathcal{C}_{2^kp}(xI_{2^kp} - \mathcal{D}_{2^kp})^{-1}\mathcal{C}^T_{2^{k}p} \end{vmatrix} \text{ is }$$ 

$$\begin{array}{@{} c @{}} 
\begin{vmatrix}
\begin{array}{ *{11}{c} } 

x-\frac{2^{k-1}p}{x} - \frac{2^{k-1}p}{x-1} & -1-\frac{2^{k-1}p}{x-1} & -1 & -1 & -1 & -1 & \cdots & -1 & -1 \\ \smallskip
-1-\frac{2^{k-1}p}{x-1} & x-\frac{2^{k-1}p}{x-1} & -1 & -1 & -1 & -1 & \cdots & -1 & -1 \\
-1 & -1 & x & -1 & -1 & -1 & \cdots & -1 & -1 \\
-1 & -1 & -1 & x & -1 & -1 & \cdots & -1 & -1 \\
-1 & -1 & -1 & -1 & x & -1 & \cdots & -1 & -1 \\
-1 & -1 & -1 & -1 & -1 & x & \cdots & -1 & -1 \\

\vdots & \vdots & \vdots & \vdots & \vdots & \vdots & \ddots & \vdots & \vdots \\

-1 & -1 & -1 & -1 & -1 & -1 & \cdots & \ x & -1 \\

-1 & -1 & \undermat{2^kp-2}{-1 & -1 & -1 & -1 & -1 & -1 &  x & } \\
\end{array}
\end{vmatrix} \\
\mathstrut
\end{array}.$$

$$ \text{ On applying successively } C_3 \rightarrow C_3 + xC_1, C_i \rightarrow C_i - C_1, \text{ where } 2 \leq i \neq 3 \leq 2^kp,$$ 
we have  \bigskip
$\begin{vmatrix}
xI_{2^kp}-\mathcal{B}_{2^kp}-\mathcal{C}_{2^kp}(xI_{2^kp} - \mathcal{D}_{2^kp})^{-1}\mathcal{C}^T_{2^{k}p} \end{vmatrix} =$
	
$$ \begin{vmatrix}

\begin{array}{ *{10}{c} }

x-\frac{2^{k-1}p}{x} - \frac{2^{k-1}p}{x-1} & -x-1+\frac{2^{k-1}p}{x} & \frac{x^3-x^2-(2^{k}p+1)x+(2^{k-1}p+1)}{x-1} \\\\
-\frac{x+2^{k-1}p-1}{x-1} & x+1 & \frac{-x^2-2^{k-1}px+1}{x-1}  \\
-1 & 0 & 0   \\
-1 & 0 & -x-1  &  & \mathcal{X}_{2^{k}p-3}\\
-1 & 0 & -x-1 \\
\vdots & \vdots & \vdots \\
-1 & 0 & -x-1   \\
\end{array}
\end{vmatrix} \\
\mathstrut,$$ where

$$\mathcal{X}_{2^{k}p-3} = \begin{bmatrix}

\begin{array}{ *{10}{c} }

 \frac{-x^3+(2^{k}p+1)x-2^{k-1}p}{x^2-x}  & \frac{-x^3+(2^{k}p+1)x-2^{k-1}p}{x^2-x} &  \cdots & \frac{-x^3+(2^{k}p+1)x-2^{k-1}p}{x^2-x} \\\\
 \frac{2^{k-1}p}{x-1} &  \frac{2^{k-1}p}{x-1} & \cdots & \frac{2^{k-1}p}{x-1}  \\
 0  & 0 & \cdots & 0  \\
 x+1 & 0 & \cdots & 0   \\
 0 & x+1 & \cdots & 0   \\
 \vdots & \vdots  & \ddots  & \vdots \\
 0  & 0 & \cdots & x+1  \\
\end{array}
\end{bmatrix} \\
\mathstrut.$$

\bigskip

On expanding along $R_3$, we have

\bigskip
\hspace{.05cm}
$$ = (-1)(-1)^{3+1} \begin{vmatrix}

\begin{array}{ *{10}{c} }
 -x-1+\frac{2^{k-1}p}{x} & \frac{x^3-x^2-(2^kp+1)x+(2^{k-1}p+1)}{x-1}  \\\\
 x+1 & \frac{-x^2-2^{k-1}px+1}{x-1} \\
 0 & -x-1 &  & \mathcal{X}_3 & \\
 0 & -x-1 \\
\vdots & \vdots \\
 0 & -x-1 \\
\end{array}
\end{vmatrix} \\
\mathstrut,$$

\smallskip

where $\mathcal{X}_3$ is the matrix obtained by deleting $R_3$ from $\mathcal{X}_{2^{k}p-3}$.

\bigskip
On applying $C_3 \rightarrow C_3 + C_2$ and expanding along $R_3$, we have
$$ = (-1
)(-1)^{3+2}(-x-1) \begin{vmatrix}

\begin{array}{ *{10}{c} }
-x-1+\frac{2^{k-1}p}{x} & \frac{x^4-2x^3-(2^kp+1)x^2+(3.2^{k-1}p+2)x-2^{k-1}p}{x^2-x} \\\\
x+1 &  \frac{-x^2-2^{k-1}px+(2^{k-1}p+1)}{x-1} \\
0  & -x-1    \\
0  & -x-1 &  \mathcal{X}'_3 \\
\vdots & \vdots  \\
0  & -x-1   \\
\end{array}
\end{vmatrix} \\
\mathstrut,$$ where $\mathcal{X}'_3$ is the matrix obtained by deleting $R_3$ from $\mathcal{X}_3.$

\bigskip
On applying $C_3 \rightarrow C_3 + C_2$ and expanding along $R_3$ successively, we have

\bigskip

$ = -(x+1)^{2^{k}p-3} \begin{vmatrix}
\begin{array}{ *{9}{c} }

 -x-1+\frac{2^{k-1}p}{x} & q(x) \\
\\
x+1 & \frac{-x^2-2^{k-1}px+[(2^{k-1}p+1)+2^{k-1}p(2^{k}p-4)]}{x-1} \\
\\

\end{array}
\end{vmatrix}, \\
\mathstrut$ where 
$q(x) =  \frac{x^4+[-2-(2^{k}p-4)]x^3-(2^kp+1)x^2+[(3.2^{k-1}p+2)+(2^{k}p-4)(2^kp+1)]x-[2^{k-1}p+2^{k-1}p(2^{k}p-4)]}{x^2-x}$

\bigskip
\begin{align*}
&= \frac{-(x+1)^{2^{k}p-3}}{x(x-1)} \Big[-x^5 + (2^kp-2)x^4 + 5.2^{k-1}px^3+(-3.2^{2k-1}p^2+2^{k+2}p+2)x^2\\
&+(1-5.2^{2k-2}p^2+3.2^{k-1}p)x+(2^{3k-2}p^3-2^{2k-2}p^2-2^kp)\Big].
\end{align*}

Therefore, 
\begin{align*}
\Delta(P(\mathcal{G}), x) &= x^{2^{k-1}p-1}(x+1)^{5.2^{k-2}p-3}(x-1)^{2^{k-2}p-1}\bigg[x^5 -(2^kp-2)x^4 - 5.2^{k-1}px^3 \\ &+(3.2^{2k-1}p^2-2^{k+2}p-2)x^2+(5.2^{2k-2}p^2-3.2^{k-1}p-1)x \\
&-(2^{3k-2}p^3-2^{2k-2}p^2-2^kp)\bigg]. \hspace{7.5cm}\hfill \Box
\end{align*}

\begin{lem}\label{l9}
	The spectral radius of $P(\mathcal{G})$ is given as:
	$$ \lambda_1(P(\mathbb{Z}_{2^{k}p})) < \lambda_1(P(\mathcal{G})) \leq \lambda_1(P(\mathbb{Z}_{2^{k}p})) + \sqrt{{2^kp}} + \frac{1+\sqrt{2^{k+1}p}}{2}.$$
\end{lem}

\noindent{\bf Proof.}
The adjacency matrix of $P(\mathcal{G})$ is $A = \mathcal{Y}+\mathcal{Z}$, where 
$\mathcal{Y} = \begin{bmatrix}
\mathcal{B}_{2^{k}p} & \mathcal{W}_{2^{k}p}  \\ \\
\mathcal{W}_{2^{k}p} & \mathcal{O}_{2^{k}p} 
\end{bmatrix}$ and $\mathcal{Z} =
\begin{bmatrix}

\mathcal{O}_{2^{k}p} & \mathcal{X}_{2^{k}p}  \\ \\
\mathcal{X}^T_{2^{k}p} & \mathcal{D}_{2^{k}p} 
\end{bmatrix}$ are Hermitian matrices. 
Here, the matrices  $\mathcal{B}_{2^{k}p}$ and $\mathcal{D}_{2^{k}p}$ are the same as in Lemma \ref{l6}.

%\hspace{-.75cm}
$$\mathcal{W}_{2^kp} =
\begin{bmatrix}
1 & 1 & 1 & \cdots & 1 \\
0 & 0 & 0 & \cdots & 0 \\
0 & 0 & 0 & \cdots & 0 \\
\vdots & \vdots & \vdots &\ddots & \vdots \\
0 & 0 & 0 & \cdots & 0
\end{bmatrix} \text{ and  } \mathcal{X}_{2^kp} =
\begin{bmatrix}
0 & 0 & 0 & 0  & \cdots & 0 \\
1 & 1 & 1 & 0  & \cdots & 0 \\
0 & 0 & 0 & 0  & \cdots & 0 \\
\vdots & \vdots & \vdots & \vdots &\ddots & \vdots \\
\undermat{2^{k-1}p}{0 & 0 & 0}  & 0  & \cdots & 0
\end{bmatrix}.$$

\bigskip
By Theorem \ref{t6.3}, we get 
\begin{equation}
\lambda_1{(A)} \leq \lambda_1(\mathcal{Y}) + \lambda_1(\mathcal{Z}).
\end{equation}

\bigskip
Now, we can write $\mathcal{Y} = \mathcal{Y}_1 + \mathcal{Y}_2$, where $\mathcal{Y}_1 =
\begin{bmatrix}

\mathcal{B}_{2^{k}p} & \mathcal{O}_{2^{k}p}  \\ \\
\mathcal{O}_{2^{k}p} & \mathcal{O}_{2^{k}p} 
\end{bmatrix}$ and $\mathcal{Y}_2 =
\begin{bmatrix}

\mathcal{O}_{2^{k}p} &  \mathcal{W}_{2^{k}p}  \\ \\
\mathcal{W}_{2^{k}p} &  \mathcal{O}_{2^{k}p} 
\end{bmatrix}.$

\bigskip
By Theorem \ref{t6.3}, we have

\begin{equation}\label{Eq4}
\lambda_1(\mathcal{Y}) \leq \lambda_1(\mathcal{Y}_1) + \lambda_1(\mathcal{Y}_2).
\end{equation}

The characteristic polynomial of $\mathcal{Y}_2$ is

$$\begin{vmatrix}
xI_{2^{k+1}p}-\mathcal{Y}_2
\end{vmatrix} =
\begin{vmatrix}

xI_{2^{k}p} & -\mathcal{W}_{2^{k}p}  \\ \\
-\mathcal{W}^T_{2^{k}p} & xI_{2^{k}p} 
\end{vmatrix}.$$

First we apply $R_1 \rightarrow xR_1$ and then $R_1 \rightarrow R_1+ R_{2^{k}p+1}$, $R_1 \rightarrow R_1 + R_{2^{k}p+2}, \ldots, R_1 \rightarrow R_1 + R_{2^{k+1}p}$.  In the resultant matrix on expanding along $R_1,$ we get

\begin{align*}
\begin{vmatrix}
xI_{2^{k+1}p}-\mathcal{Y}_2
\end{vmatrix} & = \frac{(x^2-2^{k}p)}{x} 
\begin{vmatrix}
xI_{2^{k}p-1 \times 2^{k}p-1} & \mathcal{O}_{2^{k}p-1 \times 2^{k}p}  \\
-\mathcal{W}^T_{2^{k}p \times 2^{k}p-1} & xI_{2^{k}p \times 2^{k}p}  
\end{vmatrix} \\ & = 
\frac{(x^2-2^{k}p)}{x}|xI_{2^{k}p-1 \times 2^{k}p-1}||xI_{2^{k}p \times 2^{k}p}|.
\end{align*}

Therefore, $$\begin{vmatrix}
xI_{2^{k+1}p}-\mathcal{Y}_2
\end{vmatrix} = (x^2-2^{k}p)x^{2^{k+1}p-2}.$$

\smallskip
Thus, the eigenvalues of $\mathcal{Y}_2$ are $  \pm \sqrt{2^{k}p}, \ 0^{2^{k+1}p-2}.$  
Therefore, $\lambda_1(\mathcal{Y}_2) = \sqrt{{2^{k}p}}.$

\smallskip

Note that, the adjacency matrix of $P(\mathbb{Z}_{2^kp})$ is $\mathcal{B}_{{2^k}p}$. 
Hence, $\lambda_1(\mathcal{Y}_1) = \lambda_1(P(\mathbb{Z}_{{2^k}p})).$

\smallskip
On putting the values of $\lambda_1(\mathcal{Y}_1)$ and $\lambda_1(\mathcal{Y}_2)$ in Equation (\ref{Eq4}), we get 
\begin{equation} \label{Eq5}
\lambda_1(\mathcal{Y}) \leq  \lambda_1(\mathcal{Y}_1) + \lambda_1(\mathcal{Y}_2) =  \lambda_1(P(\mathbb{Z}_{{2^k}p})) + \sqrt{{2^{k}p}}.
\end{equation}

Next, the characteristic polynomial of $\mathcal{Z}$ is 

$$\begin{vmatrix}
xI_{2^{k+1}p}-\mathcal{Z}
\end{vmatrix} =
\begin{vmatrix}

xI_{2^{k}p} & -\mathcal{X}_{2^{k}p}  \\ \\
-\mathcal{X}^T_{2^{k}p} & xI_{2^{k}p}-\mathcal{D}_{2^kp} 
\end{vmatrix}.$$

As discussed in Lemma \ref{l6} that, $xI_{2^kp}-D_{2^kp}$ is invertible. So, by Theorem \ref{t6.1},
\begin{align*}
\begin{vmatrix}
xI_{2^{k+1}p}-\mathcal{Z}
\end{vmatrix} & = \begin{vmatrix}
xI_{2^kp}-\mathcal{D}_{2^kp}
\end{vmatrix}\begin{vmatrix}
xI_{2^kp}-\mathcal{X}_{2^kp}(xI_{2^kp} - \mathcal{D}_{2^kp})^{-1}\mathcal{X}^T_{2^{k}p}
\end{vmatrix} \\ & = (x^2-1)^{2^{k-2}p}x^{{2^{k-1}}p}  \begin{vmatrix}
xI_{2^kp}-\mathcal{X}_{2^kp}(xI_{2^kp} - \mathcal{D}_{2^kp})^{-1}\mathcal{X}^T_{2^{k}p}\end{vmatrix}.
\end{align*}

Now, $xI_{2^kp}-\mathcal{X}_{2^kp}(xI_{2^kp} - \mathcal{D}_{2^kp})^{-1}\mathcal{X}^T_{2^{k}p}$ is a diagonal matrix such that its second diagonal entry is $x-\frac{2^{k-1}p}{x-1}$ and the remaining diagonal entries are $x$.

\bigskip
Therefore, $\begin{vmatrix}
xI_{2^kp}-\mathcal{X}_{2^kp}(xI_{2^kp} - \mathcal{D}_{2^kp})^{-1}\mathcal{X}^T_{2^{k}p} 
\end{vmatrix} =\frac{x^{2^{k}p-1} \left(x^2-x-2^{k-1}p \right)}{x-1}.$

Hence, $$\begin{vmatrix}
xI_{2^{k+1}p}-\mathcal{Z}
\end{vmatrix} = x^{3.2^{k-1}p-1} \left(x-1 \right)^{2^{k-2}p-1}\left(x+1 \right)^{{2^{k-2}p}}\left(x^2-x-2^{k-1}p \right).$$
Thus, the eigenvalues of $\mathcal{Z}$ are 
$$\frac{1 \pm \sqrt{1+2^{k+1}p}}{2}, \ 0^{3.2^{k-1}p-1}, \ 1^{2^{k-2}p-1}, \ (-1)^{2^{k-2}p}.$$ 
This gives that, $$ \lambda_1(\mathcal{Z}) = \frac{1 + \sqrt{1+2^{k+1}p}}{2}.$$
Also, by Theorem \ref{t6.2}, we have
$$\lambda_1(P(\mathbb{Z}_{2^{k}p})) < \lambda_1{(P(\Gamma))}.$$ 
Combining all the inequalities, we have the lemma.  \hfill $\Box$

%$$ \lambda_1(P(\mathbb{Z}_{2^{k}p})) < \lambda_1{(P(\mathcal{G}))} \leq \lambda_1(P(\mathbb{Z}_{2^{k}p}) + \sqrt{{2^{k}p}} + \frac{1 + \sqrt{1+2^{k+1}p}}{2} \hspace{3cm} $$

\begin{lem} \label{l7}
The characteristic polynomial of the Laplacian matrix $L(P(\mathcal{G}))$ of $P(\mathcal{G})$ is:

$$\Delta(L(P(\mathcal{G})), x) = x(x-2^{k+1}p)(x-3.2^{k-1}p)(x-2^kp)^{2^kp-3}(x-4)^{2^{k-2}p}(x-2)^{2^{k-2}p}
(x-1)^{2^{k-1}p}.$$ 

\end{lem}

\noindent{\bf Proof.}  The Laplacian matrix  $L(P(\mathcal{G}))$ of $P(\mathcal{G})$ is 
$\begin{bmatrix}
\mathcal{L}_{2^kp} & \mathcal{M}_{2^kp}  \\ \\
 \mathcal{M}^T_{2^kp} &  \mathcal{N}_{2^kp}
\end{bmatrix}$, where

\hspace{-.85cm}
$$\mathcal{M}_{2^kp} = \begin{array}{@{} c @{}}
\begin{bmatrix}
\begin{array}{ *{8}{c} }

\mathcal{P} & \mathcal{P} & \cdots & \mathcal{P} & \mathcal{Q} &  \mathcal{Q} & \cdots &  \mathcal{Q} \\
\mathcal{O} &  \mathcal{O} & \cdots &  \mathcal{O} &  \mathcal{O} &  \mathcal{O} & \cdots &  \mathcal{O} \\
\mathcal{O} &  \mathcal{O} & \cdots &  \mathcal{O} &  \mathcal{O} &  \mathcal{O} & \cdots &  \mathcal{O} \\
\vdots & \vdots & \vdots & \vdots & \vdots & \vdots & \ddots  & \vdots \\
\mathcal{O} &  \mathcal{O} & \cdots &  \mathcal{O} &  \mathcal{O} &  \mathcal{O} & \cdots &  \mathcal{O} \\

\undermat{2^{k-2}p}{\mathcal{O} & \mathcal{O} & \cdots & \mathcal{O}} & \undermat{2^{k-2}p}{\mathcal{O} & \mathcal{O} & \cdots & \mathcal{O}}  \\
\end{array}
\end{bmatrix} \\
\mathstrut
\end{array},   \mathcal{N}_{2^kp} = \begin{array}{@{} c @{}}
\begin{bmatrix}
\begin{array}{ *{9}{c} }

\mathcal{R} & \mathcal{O} & \mathcal{O} & \cdots & \mathcal{O} & \mathcal{O} &\cdots & \mathcal{O} \\
\mathcal{O} & \mathcal{R} & \mathcal{O} & \cdots & \mathcal{O} & \mathcal{O} & \cdots & \mathcal{O} \\
\mathcal{O} & \mathcal{O} & \mathcal{R} & \cdots & \mathcal{O} & \mathcal{O} & \cdots & \mathcal{O}\\
\vdots & \vdots & \vdots & \ddots & \vdots & \vdots & \ddots & \vdots \\
\mathcal{O} & \mathcal{O} & \mathcal{O} & \cdots & \mathcal{R} & \mathcal{O} &  \cdots & \mathcal{O}\\
\mathcal{O} & \mathcal{O} & \mathcal{O} & \cdots & \mathcal{O} & \mathcal{S} & \cdots & \mathcal{O} \\
\vdots & \vdots & \vdots & \vdots & \vdots & \vdots & \ddots & \vdots \\
\undermat{2^{k-2}p}{\mathcal{O} & \mathcal{O} & \mathcal{O} & \cdots & \mathcal{O}} & \undermat{2^{k-2}p}{\mathcal{O} & \cdots & \mathcal{S}}  \\
\end{array}
\end{bmatrix} \\
\mathstrut
\end{array},$$ $$ \mathcal{L}_{2^kp} =
\begin{bmatrix}
2^{k+1}p-1 & -1 & -1 & -1 & \cdots & -1 \\
-1 & 3.2^{k-1}p-1 & -1 & -1 & \cdots & -1 \\
-1 & -1 & 2^{k}p-1 & -1 & \cdots & -1 \\
-1 & -1 & -1 & 2^{k}p-1 & \cdots & -1 \\
\vdots & \vdots & \vdots &\vdots & \ddots &\vdots \\
-1 & -1 & -1 & -1 & \cdots & 2^{k}p-1
\end{bmatrix},$$ \  $$ \mathcal{P} =
\begin{bmatrix}
-1 & -1 \\
-1 & -1
\end{bmatrix}, \  \mathcal{Q} =
\begin{bmatrix}
-1 & -1 \\
0 & 0
\end{bmatrix},\  \mathcal{R} =
\begin{bmatrix}
3 & -1 \\
-1 & 3
\end{bmatrix}, \  \mathcal{S} =
\begin{bmatrix}
1 & 0 \\
0 & 1
\end{bmatrix}, \text{ and }\  \mathcal{O} =
\begin{bmatrix}
0 & 0 \\
0 & 0
\end{bmatrix}.$$

Therefore, the characteristic polynomial of $L(P(\mathcal{G}))$ is
$\begin{vmatrix}
xI_{2^kp}-
\mathcal{L}_{2^kp} & -\mathcal{M}_{2^kp}  \\
-\mathcal{M}^T_{2^kp} & xI_{2^kp} - \mathcal{N}_{2^kp}
\end{vmatrix}.$

\bigskip
Applying $R_1 \rightarrow (x-1)R_1$, and then $R_1 \rightarrow
R_1-R_2-R_3- \cdots -R_{2^{k+1}p}$ and expending along $R_1$, we have

$$ \frac{x(x-2^{k+1}p)}{(x-1)}
\begin{vmatrix}
\begin{array}{ *{14}{c} }
x-(3.2^{k-1}p-1) & \undermat{2^{k}p-2}{1 & 1 & \cdots & 1 &} &  \undermat{2^{k-1}p}{1 & 1 & \cdots & 1} & \undermat{2^{k-1}p}{0 & 0 & \cdots & 0}\\
1 &  &  &  &  \\
1 &  &  &  &  \\
1 &  & \mathcal{T} &  & & & & \mathcal{U} &  & &   \\
1 &  & &  &  & & &  \\

\vdots &  &  &  &  &  &  &  &  &  \\
\vdots &  &  &  &  &  &  &  &  &  \\
1 &  &  &  &  &  &  &  &  &  \\
1 &  & \mathcal{U}^T &  &  & & & \mathcal{V}  &  &  &  \\
1 &  &  &  &  &  &  &  &  & 
\end{array}
\end{vmatrix},$$

where  $ \mathcal{T}= (x-2^kp+1)I_{2^kp-2} + \mathcal{J}_{2^kp-2}, \ \mathcal{U} = \mathcal{O}_{2^kp-2 \times 2^{k}p}, \ \text{and} \  \mathcal{V} = xI_{2^kp} - \mathcal{N}_{2^kp}$.

\bigskip
Applying $R_1 \rightarrow (x-2)R_1$, and then $R_1 \rightarrow
R_1-R_2-R_3- \cdots -R_{3.2^{k-1}p-1}$ and expending along $R_1$, we have

$$ \begin{vmatrix} xI_{2^{k+1}p} - L(P(\mathcal{G}))  \end{vmatrix} = \frac{x(x-2^{k+1}p)(x-1)(x-3.2^{k-1}p)}{(x-1)(x-2)}
\begin{vmatrix}
\begin{array}{ *{14}{c} }
 \mathcal{T} &  \mathcal{U}  \\
 \mathcal{U}^T  & \mathcal{V}   \\
\end{array}
\end{vmatrix}.$$

By Theorem \ref{l6.3}, we have

$$ \begin{vmatrix} xI_{2^{k+1}p} - L(P(\mathcal{G}))  \end{vmatrix} = \frac{x(x-2^{k+1}p)(x-1)(x-3.2^{k-1}p)}{(x-1)(x-2)}
\begin{vmatrix} \mathcal{T} \end{vmatrix}
\begin{vmatrix} \mathcal{V} \end{vmatrix}.$$
Now, using Theorem \ref{l6.4}, we get
\begin{align*}
\begin{vmatrix} \mathcal{T} \end{vmatrix} &= [(x-2^kp+1) + (2^kp-2)-1][(x-2^kp+1)-1]^{(2^kp-2)-1} \\
&= (x-2)(x-2^kp)^{2^kp-3},\\
\begin{vmatrix}  \mathcal{V} \end{vmatrix} &= [(x-3)^2-1]^{2^{k-2}p}(x-1)^{2^{k-1}p}.
\end{align*}

Thus, 
$$\Delta(L(P(\mathcal{G})), x) = x(x-2^{k+1}p)(x-3.2^{k-1}p)(x-2^kp)^{2^kp-3}(x-4)^{2^{k-2}p}(x-2)^{2^{k-2}p}
(x-1)^{2^{k-1}p}.$$
Which is the required characteristics polynomial of the Laplacian matrix $L(P(\mathcal{G}))$.  \hfill$\Box$

\smallskip

The roots $\mu_1=0$, $\mu_2=1$, $\mu_3=2$, $\mu_4=4$, $\mu_5=2^k p$, $\mu_6= 3. 2^{k-1} p$, and $\mu_7= 2^{k+1} p$ of the above polynomial $\Delta(Q(P(\mathcal{G})), x)$ are called Laplacian eigenvalues. Then we have the following corollaries.

\begin{cor}
	The Laplacian spectrum of $P(\mathcal{G})$ is given below:

$$\begin{pmatrix}
	0 & 1 & 2 & 4 & 2^kp & 3.2^{k-1}p & 2^{k+1}p \\
    1 & 2^{k-1}p & 2^{k-2}p & 2^{k-2}p & 2^kp-3 & 1 & 1
\end{pmatrix}.$$

\end{cor}

\begin{cor}
The Laplacian energy $LE(P(\mathcal{G}))$ of  $P(\mathcal{G})$ is $LE(P(\mathcal{G})) = \frac{5. 2^k p - 13}{4}$.
\end{cor}

\noindent{\bf Proof.} Let $m$ and $n$ denote the number of edges and vertices in $P(\mathcal{G})$ respectively. Then we see that $m =2^{k-2} p (5+2^{k+1} p)$ and $n = 2^{k+1} p$. Now, the proof follows from the definition $LE(\mathcal{G}) = \sum_{i=1}^7|\mu_i - \frac{2m}{n}|$, where $\mu_i$, $1 \leq i \leq 7$, are Laplacian eigenvalues. \hfill$\Box$

\begin{lem} \label{l7}
	The characteristic polynomial of the signless Laplacian matrix $Q(P(\mathcal{G}))$ of $P(\mathcal{G})$ is:
	
	$\Delta(Q(P(\mathcal{G})), x) = (x-1)^{2^{k-1}p-1}(x-2)^{2^{k-2}p}(x-4)^{2^{k-2}p-1}(2^kp-x-2)^{2^{k}p-3}\bigg[-x^{5} + (11.2^{k-1}p -1)x^{4} - (5.2^{2k+1}p^2 + 5.2^{k-1}p-14)x^3 + (
 3.2^{3k+1}p^3 + 11.2^{2k+1}p^2 - 65.2^{k}p + 28)x^2 - (25.2^{3k}p^3 -135.2^{2k-1}p^2 + 9.2^{k+2}p + 8
)x + (15.2^{3k}p^3 - 31.2^{2k+1}p^2 + 5.2^{k+4}p - 32 )\bigg].$
\end{lem}

\noindent{\bf Proof.}  The signless Laplacian matrix $Q(P(\mathcal{G}))$ of $P(\mathcal{G})$ is 
\smallskip
$\begin{bmatrix}
L_{2^kp} & M_{2^kp}  \\ \\
M^T_{2^kp} & N_{2^kp}
\end{bmatrix}$, where

\hspace{-.85cm}
$$M_{2^kp} = \begin{array}{@{} c @{}}
\begin{bmatrix}
\begin{array}{ *{8}{c} }

P & P & \cdots & P & Q & Q & \cdots & Q \\
O & O & \cdots & O & O & O & \cdots & O \\
O & O & \cdots & O & O & O & \cdots & O \\
\vdots & \vdots & \vdots & \vdots & \vdots & \vdots & \ddots  & \vdots \\
O & O & \cdots & O & O & O & \cdots & O \\
\undermat{2^{k-2}p}{O & O & \cdots & O} & \undermat{2^{k-2}p}{O & O & \cdots & O}  \\
\end{array}
\end{bmatrix} \\
\mathstrut
\end{array},~ N_{2^kp} = \begin{array}{@{} c @{}}
\begin{bmatrix}
\begin{array}{ *{9}{c} }

R & O & O & \cdots & O & O &\cdots & O \\
O & R & O & \cdots & O & O & \cdots & O\\
O & O & R & \cdots & O & O & \cdots & O\\
\vdots & \vdots & \vdots & \ddots & \vdots & \vdots & \ddots & \vdots \\
O & O & O & \cdots & R & O &  \cdots & O\\
O & O & O & \cdots & O & S & \cdots & O \\
\vdots & \vdots & \vdots & \vdots & \vdots & \vdots & \ddots & \vdots \\
\undermat{2^{k-2}p}{O & O & O & \cdots & O} & \undermat{2^{k-2}p}{O & \cdots & S}  \\
\end{array}
\end{bmatrix} \\
\mathstrut
\end{array},$$

 $$L_{2^kp} =
\begin{bmatrix}
2^{k+1}p-1 & 1 & 1 & 1 & \cdots & 1 \\
1 & 3.2^{k-1}p-1 & 1 & 1 & \cdots & 1 \\
1 & 1 & 2^{k}p-1 & 1 & \cdots & 1 \\
1 & 1 & 1 & 2^{k}p-1 & \cdots & 1 \\
\vdots & \vdots & \vdots &\vdots & \ddots &\vdots \\
1 & 1 & 1 & 1 & \cdots & 2^{k}p-1
\end{bmatrix},  \ P =
\begin{bmatrix}
1 & 1 \\
1 & 1
\end{bmatrix},$$ $$ \ Q =
\begin{bmatrix}
1 & 1 \\
0 & 0
\end{bmatrix},\ R =
\begin{bmatrix}
3 & 1 \\
1 & 3
\end{bmatrix}, \ S =
\begin{bmatrix}
1 & 0 \\
0 & 1
\end{bmatrix}, \ O =
\begin{bmatrix}
0 & 0 \\
0 & 0
\end{bmatrix}.$$
\\
Therefore, we have $ \Delta(Q(P(\mathcal{G})), x)=
\begin{vmatrix}
xI_{2^kp}-
L_{2^kp} & - M_{2^kp}  \\
-M^T_{2^kp} & xI_{2^kp} - N_{2^kp}
\end{vmatrix} .$ Now, by Proposition \ref{t6}, we have

 $$(xI_{2^kp} - N_{2^kp})^{-1} = \begin{array}{@{} c @{}}

\begin{bmatrix}
\begin{array}{ *{9}{c} }

\frac{x-3}{(x-2)(x-4)} & \frac{1}{(x-2)(x-4)} & \cdots & 0 &  0 & 0 & \cdots & 0  \\
\frac{1}{(x-2)(x-4)} & \frac{x-3}{(x-2)(x-4)} & \cdots & 0 & 0 & 0 & \cdots & 0 \\
\vdots & \vdots & \ddots & \vdots & \vdots & \vdots & \ddots & \vdots  \\
0 & 0 & \cdots & \frac{x-3}{(x-2)(x-4)} & \frac{1}{(x-2)(x-4)} & 0 & \cdots & 0  \\
0 & 0 & \cdots & \frac{1}{(x-2)(x-4)} & \frac{x-3}{(x-2)(x-4)} & 0 & \cdots & 0  \\
0 & 0 & \cdots & 0 & 0 & \frac{1}{x-1} & \cdots & 0 \\
\vdots & \vdots & \ddots & \vdots &  \vdots & \vdots & \ddots & \vdots  \\
0 & 0 & \cdots & 0 & 0 & \undermat{2^{k-1}p}{0 & \cdots & \frac{1}{x-1}}
\end{array}
\end{bmatrix} \\
\mathstrut
\end{array}.$$

So, by Theorem \ref{t6.1},
$$\Delta(Q(P(\mathcal{G})), x) = \begin{vmatrix}
xI_{2^kp}- N_{2^kp}
\end{vmatrix}\begin{vmatrix}
xI_{2^kp}-
L_{2^kp}-M
_{2^kp}(xI_{2^kp} - N_{2^kp})^{-1}M^T_{2^{k}p}
\end{vmatrix}.$$
Now, 
$$M_{2^kp}(xI_{2^kp} - N_{2^kp})^{-1}M^T_{2^kp} = \begin{array}{@{} c @{}}

\begin{bmatrix}
\begin{array}{ *{9}{c} }

\frac{2^{k-1}p}{x-1} + \frac{2^{k-1}p}{x-4} & \frac{2^{k-1}p}{x-4} & 0 & \cdots & 0 & 0  \\ \\
\frac{2^{k-1}p}{x-4} & \frac{2^{k-1}p}{x-4} & 0 & \cdots & 0 & 0 \\
0 & 0 & 0 & \cdots & 0 & 0 \\
\vdots & \vdots & \vdots & \vdots & \vdots & \vdots \\
0 & 0 & 0 & \cdots & 0 & 0 \\
0 & 0 & \undermat{2^kp-2}{ 0 & \cdots & 0 & 0} \\
\end{array}
\end{bmatrix} \\
\mathstrut
\end{array}.$$

Therefore, $\begin{vmatrix}
xI_{2^kp}-
L_{2^kp}-M
_{2^kp}(xI_{2^kp} - N_{2^kp})^{-1}M^T_{2^{k}p} \end{vmatrix}$ is
$$\begin{array}{@{} c @{}} 
\begin{vmatrix}
\begin{array}{ *{11}{c} } 

a & \frac{-2^{k-1}p-x+4}{x-4} & -1 &  \cdots & -1 \\ \smallskip
\frac{-2^{k-1}p-x+4}{x-4} & \frac{-2^{k-1}p+(x-4)(-3.2^{k-1}p+x+1)}{x-4} & -1 &  \cdots & -1  \\
-1 & -1 & x-2^kp+1 &  \cdots & -1   \\
\vdots & \vdots & \vdots & \ddots & \vdots   \\
-1 & -1 & -1 & \cdots & x-2^kp+1 \\

\end{array}
\end{vmatrix} \\
\mathstrut
\end{array},$$ where $a = 
x+1-2^{k+1}p-(\frac{2^{k-1}p}{x-1} + \frac{2^{k-1}p}{x-4})$.

\bigskip
On applying successively $C_3 \rightarrow C_3 + (x-2^kp+1)C_1$, $C_i \rightarrow C_i - C_1$, where $2 \leq i \neq 3 \leq 2^kp$, we have  \bigskip
$\begin{vmatrix}
xI_{2^kp}-
L_{2^kp}-M
_{2^kp}(xI_{2^kp} - N_{2^kp})^{-1}M^T_{2^{k}p} \end{vmatrix} =$

$$ \begin{vmatrix}

\begin{array}{ *{10}{c} }

a &  \frac{2^{k+1}px-3.2^{k-1}p-x^2-x+2}{x-1} & 
-1 + (x-2^{k}p+1)a \\\\
\frac{-2^{k-1}p-x+4}{x-4} & -3.2^{k-1}p+x+2 & \frac{-x-(x-2^kp+1)(2^{k-1}p+x-4)+4}{x-4}  \\ \\
-1 & 0 & 0   \\
-1 & 0 & 2^{k}p-x-2  &  & X \\
-1 & 0 & 2^{k}p-x-2 \\
\vdots & \vdots & \vdots \\
-1 & 0 & 2^{k}p-x-2 \\
\end{array}
\end{vmatrix} \\
\mathstrut,$$ where $a = 
x+1-2^{k+1}p-\frac{2^{k-1}p}{x-1} - \frac{2^{k-1}p}{x-4} $ and

$$X = \begin{bmatrix}

\begin{array}{ *{10}{c} }

-a-1  & 
-a-1  &  \cdots & 
-a-1  \\\\
\frac{2^{k-1}p}{x-4} &  \frac{2^{k-1}p}{x-4} & \cdots & \frac{2^{k-1}p}{x-4}  \\
0  & 0 & \cdots & 0  \\
x-2^{k}p+2 & 0 & \cdots & 0   \\
0 & x-2^{k}p+2 & \cdots & 0   \\
\vdots & \vdots  & \ddots  & \vdots \\
0  & 0 & \cdots & x-2^{k}p+2  \\
\end{array}
\end{bmatrix} \\
\mathstrut.$$

\bigskip
 On expanding along $R_3$, we have
$$ = (-1)(-1)^{3+1} \begin{vmatrix}

\begin{array}{ *{10}{c} }

 \frac{2^{k+1}px-3.2^{k-1}p-x^2-x+2}{x-1} & 
-1 + (x-2^{k}p+1)a \\\\
 -3.2^{k-1}p+x+2 & \frac{-x-(x-2^kp+1)(2^{k-1}p+x-4)+4}{x-4}
  & &  & X_3  \\ \\
0 & 2^{k}p-x-2 \\
\vdots & \vdots \\
0 & 2^{k}p-x-2 \\
\end{array}
\end{vmatrix} \\
\mathstrut,$$ where the matrix $X_3$ is obtained by deleting $R_3$ from $X$.

\smallskip
On applying $C_3 \rightarrow C_3 + C_2$ and expanding along $R_3$, we have

$$ = 
(2^kp-x-2) \begin{vmatrix}

\begin{array}{ *{10}{c} }
\frac{2^{k+1}px-3.2^{k-1}p-x^2-x+2}{x-1} & -1+(x-2^kp+1)a + (-a-1) \\\\
-3.2^{k-1}p+x+2  &  \frac{-x-(x-2^kp+1)(2^{k-1}p+x-4)+4}{x-4} + \frac{2^{k-1}p}{x-4} \\
0  & 2^kp-x-2    \\
0  & 2^kp-x-2 & & X'_3 \\
\vdots & \vdots  \\
0  & 2^kp-x-2   \\
\end{array}
\end{vmatrix} \\
\mathstrut,$$ where $X'_3$ is the matrix obtained by deleting $R_3$ from $X_3$.

\smallskip
On applying $C_3 \rightarrow C_3 + C_2$ and expanding along $R_3$ successively, we have

\bigskip

\hspace{-.5cm}
$ = (2^kp-x-2)^{2^{k}p-3} \begin{vmatrix}
\begin{array}{ *{9}{c} }

 \frac{2^{k+1}px-3.2^{k-1}p-x^2-x+2}{x-1} & -1+(x-2^kp+1)a + (2^kp-3)(-a-1) \\
\\
-3.2^{k-1}p + x+2 & \frac{-x-(x-2^kp+1)(2^{k-1}p+x-4)+4}{x-4} + \frac{(2^kp-3)2^{k-1}p}{x-4}\\
\\

\end{array}
\end{vmatrix}, \\
\mathstrut$

$ = \frac{(2^kp-x-2)^{2^{k}p-3}}{(x-1)(x-4)}\bigg[-x^{5} + (11.2^{k-1}p -1)x^{4} - (5.2^{2k+1}p^2 + 5.2^{k-1}p-14)x^3 + (
 3.2^{3k+1}p^3 + 11.2^{2k+1}p^2 - 65.2^{k}p + 28)x^2 - (25.2^{3k}p^3 -135.2^{2k-1}p^2 + 9.2^{k+2}p + 8
)x + (15.2^{3k}p^3 - 31.2^{2k+1}p^2 + 5.2^{k+4}p - 32 )\bigg].$

\bigskip

As,  $\begin{vmatrix}
xI_{2^kp}- N_{2^kp}
\end{vmatrix}$ is a block diagonal matrix, so its determinant is the product of the determinant of its block. This gives that, $\begin{vmatrix}
xI_{2^kp}- N_{2^kp}
\end{vmatrix} = (x-1)^{2^{k-1}p}(x-2)^{2^{k-2}p}(x-4)^{2^{k-2}p}.$

\bigskip
Therefore, 

$\Delta(Q(P(\mathcal{G})), x) = (x-1)^{2^{k-1}p-1}(x-2)^{2^{k-2}p}(x-4)^{2^{k-2}p-1}(2^kp-x-2)^{2^{k}p-3}\bigg[-x^{5} + (11.2^{k-1}p -1)x^{4} - (5.2^{2k+1}p^2 + 5.2^{k-1}p-14)x^3 + (
 3.2^{3k+1}p^3 + 11.2^{2k+1}p^2 - 65.2^{k}p + 28)x^2 - (25.2^{3k}p^3 -135.2^{2k-1}p^2 + 9.2^{k+2}p + 8
)x + (15.2^{3k}p^3 - 31.2^{2k+1}p^2 + 5.2^{k+4}p - 32 )\bigg].$ \hfill $\Box$

%\section{Conclusion and future work} \label{Sec6}The present articles studied the spectral properties of power graphs of split metacyclic groups. We derived the characteristic polynomials and spectral radius. We further discussed the Laplacian energies of such graphs. The theory of graph spectra is a prominent area.  We studied several papers on spectral graph theory to understand its algebraic structure. The adjacency matrix can also be used to help with the spectral embedding of graphs in the plane. Eigenvalues of the Laplacian matrix support spectral clustering methods in machine learning. Moreover, it is included in load-balancing methods by computer scientists. Algebraic graph theory may be used to construct and investigate network topologies. Typically, supercomputer processors are put together using Cayley graphs with a high degree of symmetry.

%Nevertheless, certain eigenvalues of such graphs remain a mystery in respect of the complicated eigenvalues of the quotient matrix. In general, the eigenvalues of such graphs for other finite groups have not yet been analyzed, and extreme characterizations in terms of a variety of spectral invariants remain open and may be explained in the coming decades.

\section{Acknowledgments}The first author is thankful to the Ministry of Human Resource Development (MHRD) New Delhi, India for financial support.

\end{document}